\documentclass[11pt]{amsart}
\usepackage{geometry}                % See geometry.pdf to learn the layout options. There are lots.
\usepackage{graphicx}

\usepackage{graphics}
\usepackage{amssymb}
\usepackage{epstopdf}
\title{Self-intersections are empirically Gaussian}
\author{Moira Chas}
\begin{document}

\maketitle

\begin{abstract} In an orientable surface with boundary, free homotopy classes of curves on surfaces are in one to one correspondence with cyclic reduced words in a set of standard generators of the fundamental group. The combinatorial length of a class is the number of letters of the corresponding word. 
 The self-intersection of a free homotopy class (that is, the smallest number of self-crossings of a representative of a class) can be computed in terms of the word. 
For each of the free homotopy classes of length twenty on the punctured torus, we compute its self-intersection number and make a histogram of how many have self-intersection 0, 1, 2..... The histogram is essentially  Gaussian.
\end{abstract}
By choosing a set of free generators on the fundamental group of the punctured torus $T$, one can establish a one-to-one correspondence between the set of free homotopy classes of oriented curves on $T$ and the conjugacy classes of the free group of rank two. The \emph{combinatorial length of a conjugacy class} is the minimal number of letters (generators or inverses) required for its description.
\begin{figure}[htp]\centering
\includegraphics[width=10cm]{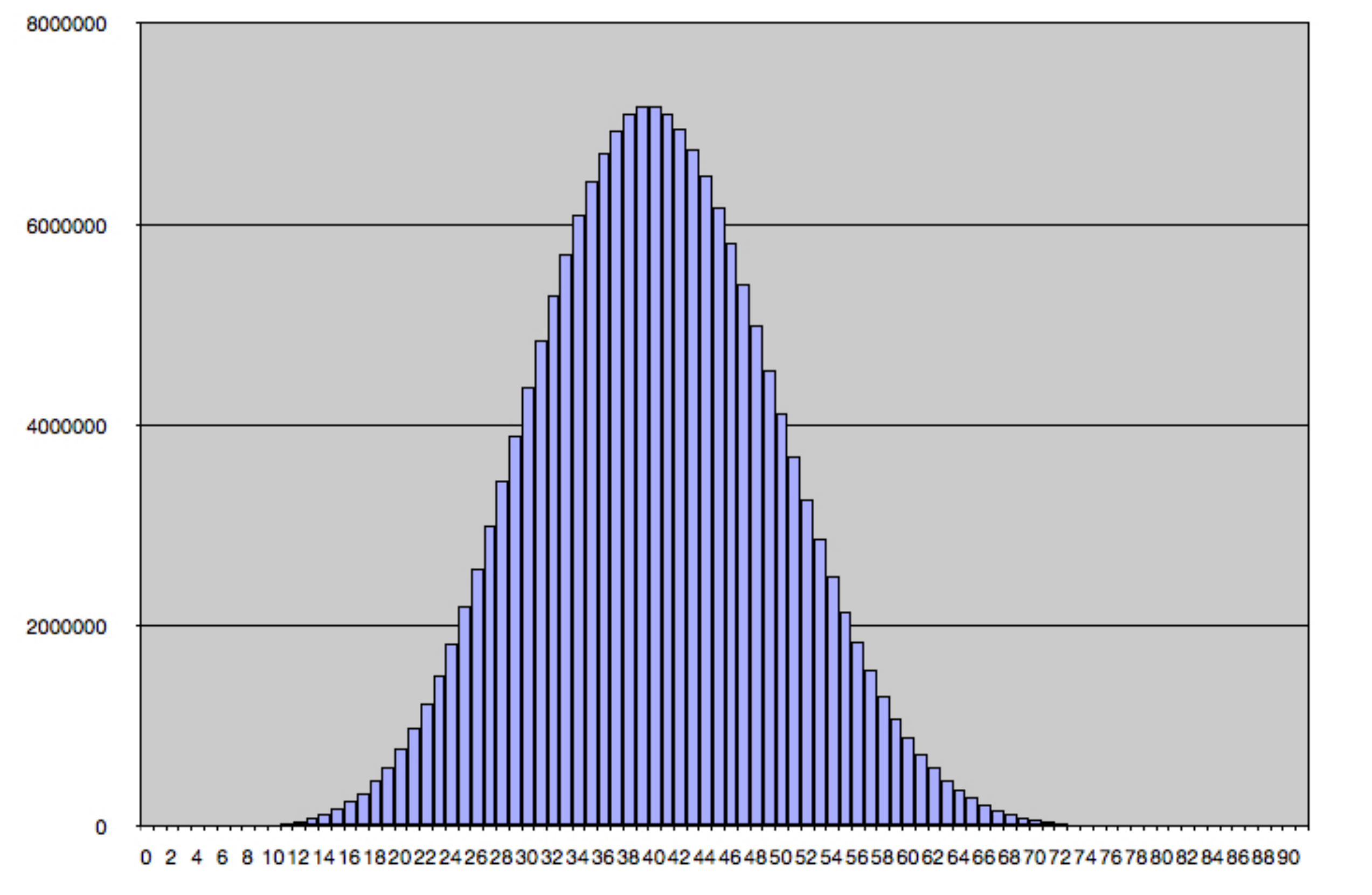} 
\caption{The histogram showing the distribution of self-intersection
numbers over  all hundred and fifty million  non-power classes of combinatorial length 20 in the punctured torus. The horizontal
coordinate shows the  self-intersection count $k$; the vertical
coordinate shows the  number of cyclic reduced words for which the
self-intersection number is  $k$. }\label{histogram}    
\end{figure}

The \emph{self-intersection of a free homotopy class of curves} (that is, the smallest number of times in which a representative crosses itself, counted by multiplicity) can be computed using the algorithms described in  \cite{cohen-lustig} or \cite{Chas:2004}. 
We programmed these algorithms using Java \cite{Chas:code} and performed this computation for each non-power class of combinatorial length twenty. Then we organized the output in the histogram of Figure \ref{histogram} portraying proportion of classes of with each possible self-intersection number. The histogram appears clearly normal. Based on these findings there is a sequel to this note \cite{Chas-Lalley} that in any orientable surface with boundary and negative Euler characteristic, the distribution of self-intersection sampling classes by their combinatorial length approaches a Gaussian as the combinatorial length goes to infinity.

\bibliographystyle{plain}

\bibliography{mainbib}

\end{document}